\documentclass[11pt, oneside]{article}

\usepackage{amsmath, amssymb, amsthm}
\usepackage{xcolor}
\usepackage{geometry, enumerate}
\usepackage[font=small,labelfont=bf, width=.85\textwidth]{caption}
\usepackage{graphicx}
\usepackage{comment}
\usepackage{verbatim}
\usepackage{mdframed}
\usepackage{enumitem}
\usepackage{capt-of}
\usepackage{ifthen}
\usepackage{epsfig}
\usepackage{tcolorbox}
\usepackage[font=small,labelfont=bf, width=.7\textwidth]{caption}

\newtheorem{theorem}{Theorem}

\newtheorem{lemma}[theorem]{Lemma}

\newtheorem{corollary}[theorem]{Corollary}
\newtheorem{conjecture}[theorem]{Conjecture}

\newtheorem{claim}{Claim}

\newtheorem{algorithm}{Algorithm}

\newenvironment{proofc}{\begin{proof}[Proof of Claim]}{\end{proof}}

\begin{document}

\title{On graphs with chromatic number and maximum degree both equal to nine}
\author{
Rachel Galindo\footnote{Auburn University, Department of Mathematics and Statistics, Auburn U.S.A.
  Email: {\tt rcg0036@auburn.edu}.}
\qquad
Jessica McDonald\footnote{Auburn University, Department of Mathematics and Statistics, Auburn U.S.A.
  Email: {\tt mcdonald@auburn.edu}.   
	Supported in part by Simons Foundation Grant \#845698  }
}

\date{}

\maketitle

\begin{abstract}  
An equivalent version of the Borodin-Kostochka Conjecture, due to Cranston and Rabern, says that any graph with $\chi = \Delta = 9$ contains $K_3 \lor E_6$ as a subgraph. Here we prove several results in support of this conjecture, where vertex-criticality and forbidden substructure conditions get us either close or all the way to containing $K_3 \lor E_6$.
\end{abstract}

\vspace*{.3in}

\section{Introduction}

In this paper all graphs are simple. For terms not defined here, we follow \cite{West}.

Given any graph $G$, it is easy to show that the chromatic number of $G$, denoted $\chi=\chi(G)$, is bounded above by $\Delta+1$, where $\Delta=\Delta(G)$ is the maximum degree of $G$.  A cornerstone of chromatic graph theory is Brooks' Theorem, which says that if $\Delta \geq 3$, then $\chi=\Delta+1$ implies that $G$ contains $K_{\Delta+1}$ as a subgraph. Borodin and Kostochka \cite{BK} conjectured the following extension of this in 1977:

\begin{conjecture}[Borodin and Kostochka \cite{BK}]
\label{BK}
     Let $G$ be a graph with $\Delta \geq 9$. Then $\chi =\Delta$ implies that $G$ contains $K_{\Delta}$ as a subgraph.
\end{conjecture}

The $\Delta\geq 9$ in Conjecture \ref{BK} is necessary: the graph obtained by blowing up $C_5$ with triangles -- that is, replacing each vertex in $C_5$ with a triangle and replacing each edge with the complete bipartite graph $K_{3,3}$ -- has $\chi=\Delta=8$ but does not contain $K_8$. Conjecture  \ref{BK} is known to be true when $\Delta>10^{14}$ (Reed \cite{Reed1999}).  It is also known to be true when induced copies of certain small graphs are forbidden, in particular: the claw $K_{1,3}$ (Cranston and Rabern \cite{CR-ClawFree}); $P_6$ (Wu and Wu \cite{WuWu}), and $P_5$ and $C_4$ (Gupta and Pradhan \cite{P5C4-free-graphs}, with an improvement by Cranston, Lafayette and Rabern \cite{(P5gemfree)}). There is also an approximate version that is known, namely the statement is true if $\Delta\geq 13$ and $K_{\Delta}$ is replaced with $K_{\Delta-3}$ (Cranston and Rabern \cite{CranRab}); this result builds on prior work by Kostochka \cite{Kostochka1980} and Mozhan \cite{Mozhan}.

In order to make progress on Conjecture \ref{BK} many authors assume that $\Delta$ is larger than its minimum possible value of nine (e.g. \cite{CranRab}, \cite{Kostochka1980}, \cite{Mozhan} mentioned above, see also Haxell and MacDonald \cite{HaxellMac}) -- one reason for this will become apparent soon when we discuss \emph{Mozhan partitions}, one of the main techniques in this area. However  Catlin \cite{Catlin} and Kostochka \cite{Kostochka1980} both independently proved that it suffices to prove Conjecture \ref{BK} for $\Delta=9$ only. Hence the topic of this paper: here we focus exclusively on graphs with chromatic number and maximum degree both equal to nine.  When $G$ is a graph with $\chi=\Delta=9$ then $G$ is known to contain $K_5$ (Borodin and Kostochka \cite{BK}); such a $G$ is also known to contain a $K_4$ in which all the vertices of the $K_4$ have degree $\Delta$ (Cranston and Rabern \cite{CranRab}). Cranston and Rabern have shown that it is not necessary to find a full $K_9$ in such graphs however, and in fact the following 
conjecture is completely equivalent to Conjecture \ref{BK}. Note that here and in what follows, given graphs $G_1, G_2$, by $G_1 \lor G_2$ we mean the graph obtained by taking disjoint copies of $G_1, G_2$ and then adding all possible edges between $G_1$ and $G_2$. By $E_t$ we mean the edgeless graph on $t$ vertices.

\begin{conjecture}[Cranston and Rabern \cite{CR}]
    \label{CR conjecture}
    Any graph with $\chi = \Delta = 9$ contains $K_3 \lor E_6$ as a subgraph.
\end{conjecture}

We have now seen that in order to prove Conjecture \ref{BK} we may restrict ourselves to graphs with   $\chi = \Delta = 9$ and we need not find a full $K_9$, but rather $K_3 \lor E_6$ suffices --- and there is yet one more significant restriction we can make. It is easy to show, using Brooks' Theorem, that it suffices to prove Conjecture \ref{BK}  (and Conjecture \ref{CR conjecture}) for \emph{vertex-critical} graphs, that is, for graphs $G$ where $\chi(G-v)<\chi(G)$ for all $v\in V(G)$.  Of course, thinking about Conjecture \ref{BK}, if $G$ is vertex-critical with $\chi=\Delta=9$ and it contains $K_9$ as a subgraph, then since $\chi(K_9)=9$ we get $G=K_9$, which is a contradiction since $\Delta(K_9)=8$. So another completely equivalent form of Conjecture \ref{BK} (and Conjecture \ref{CR conjecture}) is that there \emph{are no} vertex-critical graphs with $\chi=\Delta=9$. In this paper however, our target form of the Borodin-Kostochka Conjecture will simply be Conjecture \ref{CR conjecture} with the vertex-criticality assumption added. Our results will all have the form: if $G$ is a vertex-critical graph with $\chi = \Delta = 9$ and \emph{does not contain some special forbidden substructures}, then $G$ contains $K_3 \lor E_6$ (or is close to containing $K_3 \lor E_6$). 

\begin{figure}[htb]
    \centering
    \includegraphics[height=4cm]{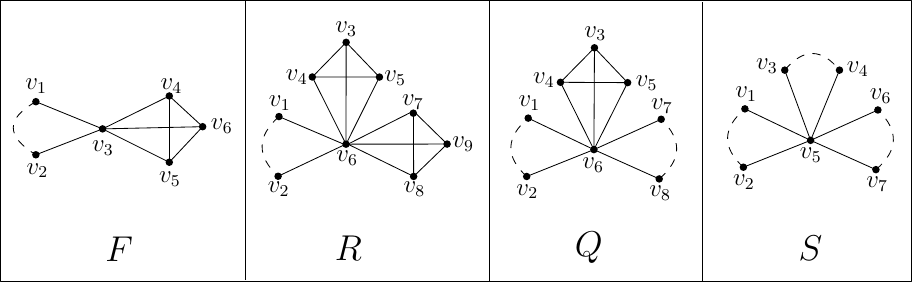}
    \caption{The forbidden substructures $F$, $R$, $Q$, and $S$.
C    Note that a path with endpoints joined by a dotted line indicates that the path is part of an induced odd cycle in the larger graph $G$. Some vertices are required to have degree either $\Delta$ or $\Delta-1$ in $G$, as listed in Table 1. }
    \label{forbidden graphs}
  \end{figure}

The forbidden substructures we consider in this paper all have one of four basic forms -- $F,R,Q$ or $S$, as depicted in Figure \ref{forbidden graphs}. In each of these figures, a path of length two with endpoints joined by a dotted line indicates that the path is part of an induced odd cycle in the larger graph $G$. The rest of the structure need not be an induced subgraph in $G$, but we do specify the degree of certain vertices to be either $\Delta$ or $\Delta-1$ in $G$. These restrictions are given in Table \ref{Table 1}; vertices not listed in the table can have any degree in $G$. As an example, a graph $G$ contains $R_2$ if it contains $R$ as a subgraph, where: $(v_1, v_6, v_2)$ happens to be part of an induced odd cycle in $G$; $v_1$ and $v_2$ can have any degree in $G$; $v_6$ has degree $\Delta -1$ in $G$, and; all other vertices in $R_2$ have degree $\Delta$ in $G$. It may be helpful to keep in mind that even vertices allowed to have any degree in $G$ only have two choices: it is easy to show that in any vertex-critical graph $G$ with $\chi=\Delta$, all vertices must have degrees $\Delta$ or $\Delta-1$.

\begin{table}
\begin{center}
    \begin{tabular}{||c|c|c||}
    \hline
        $\textcolor{white}{Graph}$ & $\Delta$-Vertices & ($\Delta-1$)-Vertices \\
        \hline 
         $F$ & $v_1, \dots, v_6$ &  \textcolor{white}{$\emptyset$} \\
         \hline
         $R_1$ & $v_1, \dots, v_5, v_7, v_8, v_9$ & $v_6$ \\
         \hline 
         $R_2$ & $v_3, v_4, v_5, v_7, v_8, v_9$ & $v_6$ \\
         \hline
         $R_3$ & $v_3, v_4, v_5, v_7, v_8, v_9$ & $v_1, v_6$ \\
         \hline 
         $R_4$ & $v_3, v_4, v_5$ & $v_1, v_6$ \\
         \hline 
         $R_5$ & $v_3, v_4, v_5$ & $v_6, v_7$ \\
         \hline
         $R_6$ & $v_3, v_4, v_5$ & $v_6$ \\
         \hline 
         $R_7$ & $v_1, v_2$ & $v_6, v_7$ \\
         \hline
         $R_8$ & $v_3, v_4$ & $v_1, v_5, v_6$ \\
         \hline
         $R_9$ & $v_3, v_4$ & $ v_5,v_6,v_7$ \\
         \hline 
         $R_{10}$ & $v_1, \dots, v_5$ & $v_6, v_7$ \\
         \hline 
         $Q$ & $v_1, \dots, v_5$ & $v_6$ \\
         \hline 
         $S$ & $v_3, v_4$ & $v_5, v_6$ \\
         \hline 
         
    \end{tabular}
    \caption{Vertices in the $F, R, Q, S$ forbidden substructures that are required to have degree $\Delta$ or $\Delta -1$ in $G$. Vertices not listed can have any degree in $G$.}
    \label{Table 1}
\end{center}
\end{table}

By forbidding certain groups of these structures, we can get that our target $G$ contains all of $K_3 \lor E_6$.

\begin{theorem} \label{fullstrong}
    Let $G$ be a vertex-critical graph with $\chi= \Delta = 9$ which does not contain any of $F, R_1, R_4, R_5, R_7, R_8, R_9$. Then $G$ contains $K_3 \lor E_6$ as a subgraph.
\end{theorem}

The graphs $R_4, R_5$, $R_7-R_9$ all contain $S$, since each contains three triangles intersecting at a vertex of degree $\Delta-1$ (with one of the triangles containing another vertex of degree $\Delta-1$ and another of the triangles containing two vertices of degree $\Delta$). So Theorem \ref{fullstrong} also yields the following corollary.

\begin{corollary}\label{fullcor} 
    Let $G$ be a vertex-critical graph with $\chi= \Delta = 9$ which does not contain any of $F, R_1, S$. Then $G$ contains $K_3 \lor E_6$ as a subgraph.
\end{corollary}

If we are okay with some edges of $K_3 \lor E_6$ being missing, then we can get get more relaxed forbidden subgraph conditions. In particular, we get the following, where 
by $H^{-k}$ we mean a graph obtained from $H$ by removing $k$ edges.

\begin{theorem}\label{minus4strong}
    Let $G$ be a vertex-critical graph with $\chi= \Delta = 9$ which does not contain any of $F, R_2, R_5$. Then $G$ contains $(K_3 \lor E_6)^{-4}$ as a subgraph.
\end{theorem}

\begin{theorem}\label{minus6strong}
    Let $G$ be a vertex-critical graph with $\chi= \Delta = 9$ which does not contain any of $F, R_1, R_3, R_{10}$. Then $G$ contains $(K_3 \lor E_6)^{-6}$ as a subgraph.
\end{theorem}

Theorems \ref{minus4strong} and \ref{minus6strong} also yield the following corollaries (since both $R_2$ and $R_5$ contain $R_6$, and  since $R_1, R_3, R_{10}$ all contain $Q$). 

\begin{corollary}\label{minus4cor}
    Let $G$ be a vertex-critical graph with $\chi= \Delta = 9$ which does not contain any of $F, R_6$. Then $G$ contains $(K_3 \lor E_6)^{-4}$ as a subgraph.
\end{corollary}

\begin{corollary}\label{minus6cor2}
    Let $G$ be a vertex-critical graph with $\chi= \Delta = 9$ which does not contain any of $F, Q$. Then $G$ contains $(K_3 \lor E_6)^{-6}$ as a subgraph.
\end{corollary}

All of Theorems \ref{fullstrong}, \ref{minus4strong} and \ref{minus6strong} are proved using the technique of \emph{Mozhan partitions}, first developed by Mozhan \cite{Mozhan} in the early 1980's. Section 2 of this paper will provide the necessary background on this method. The proofs of Theorems \ref{fullstrong} and \ref{minus4strong} are given in Section 3; Section 4 contains the proof of Theorem \ref{minus6strong}. 

Observe that none of the forbidden substructures in our above theorems are necessarily contained within either a $K_5$, or a $K_4$ having all vertices of degree $\Delta$ -- which is important, since we know those structures \emph{must} be present (in any graph with $\chi=\Delta=9$) by the above-mentioned work of Borodin and Kostochka \cite{BK} and Cranston and Rabern \cite{CranRab}, respectively. It is not obvious that a graph with $\chi= \Delta = 9$ should be \emph{able} to forbid all of the substructures listed in our results above: but it can. In the final section of this paper, Section 5, we present an infinite family of graphs with $\chi=\Delta=9$ which contain none of the substructures listed in Table 1 (or listed in our above results). Of course, our infinite family of graphs are not vertex-critical: as discussed above, finding a vertex-critical graph with $\chi=\Delta=9$ would amount to disproving Conjecture \ref{BK}. However, our forbidden substructure assumptions do not seem to get in the way of vertex-criticality. In particular, Kostochka and Yancey \cite{kostochkayancey} have established a lower bound on the average degree for vertex-critical graphs; this bound implies that a vertex-critical graph with $\chi=9$ must have average degree at least $8.75$. All the members of our infinite family have average degree higher than this, even while they forbid all our listed substructures.

\section{Mozhan Partitions}

We now introduce the partition structure which will be the main tool in our proofs. This technique was first developed by Mozhan \cite{Mozhan}, and we follow the notation used by Kostochka, Rabern, and Stiebitz in \cite{DefsPaper}, with a few minor simplifications.

Let $G$ be a vertex-critical graph with $\chi= 1+ t_1 + \dots + t_p \geq 4$, for some positive integers $t_1, \ldots, t_p, p$. Then a \emph{$(t_1, \dots, t_p)$-partition} of $G$ is a sequence $(v, X_1, \dots , X_p)$ where $v \in V(G)$ is called the \emph{special vertex}, $X_1, \dots, X_p$ is a partition of $V(G)\setminus \{v\}$, and $\chi(G[X_i]) = t_i$ for $i = 1, \dots, p$. Such a partition is easily obtained from any $\chi$-colouring of $G$ where one colour class has size one (always possible since $G$ is vertex-critical), by simply taking $v$ to be this special colour class of size one, and then grouping $t_i$ of the colour classes together to form $X_i$, for all $i$. 
This way of making a $(t_1, \dots, t_p)$-partition of $G$ means that the only edges in $G[X_i]$ are between the $t_i$ colour classes of $X_i$. A $(t_1, \dots, t_p)$-partition of $G$ is said to be \emph{optimal} if 
$|E[X_1]| + \dots + |E[X_p]|$ is minimum over all $(t_1, \dots, t_p)$-partitions of $G$, where $E[X] = E(G$[$X$]). 

Given any vertex-critical graph $G$ with $\chi\geq 1$, an old result of Dirac \cite{Dirac} tells us that $G$ is $(\chi-1)$-edge-connected, and hence that the minimum degree in $G$ is at least $\chi-1$. We call a vertex in $G$ \emph{high} if its degree is $\geq \chi$ in $G$, and \emph{low} if its degree is  $< \chi$. Note that in the case $\chi=\Delta$ we are interested in, our high vertices will all have degree $\Delta$ and our low vertices will all have degree $\Delta-1$. An optimal $(t_1, \dots, t_p)$-partition of $G$ is said to be \emph{proper} if the special vertex $v$ is low. 

The following result is the cornerstone of the method of \emph{Mozhan partitions} and is fundamental to our work in this paper; a short proof, which uses Brooks' Theorem, is written in \cite{DefsPaper}. Here and in what follows, for a vertex $a \in V(G)$ and a vertex set $X \subseteq V(G)$, we let $N_X(a) = N(a)\cap X$ and we let $d_X(a) = |N_X(a)|.$ If $(v, X_1, \dots , X_p)$ is a $(t_1, \dots, t_p)$-partition of $G$, then by $X_i^*$ we mean the component of $G[\{v\}\cup X_i]$ containing $v$, for any $1\leq i\leq p$.

\begin{lemma}[Mozhan \cite{Mozhan}]
    \label{main mozhan lemma}
    Let $G$ be a vertex-critical graph and let $(v, X_1, \dots , X_p)$ be an optimal $(t_1, \dots, t_p)$-partition of $G$. Then the following statements hold:
    \begin{enumerate}
    \item $\chi(X_i^*)=t_i+1$ and $d_{X_i}(v)\geq t_i$ for all $i$. 
        \item If $v$ is a low vertex of $G$, then $d_{X_i}(v) = t_i$ for all $i$. 
        \item If $d_{X_i}(v) = t_i$ for some $i$, then either $X_i^* = K_{t_i+1}$ or $t_i = 2$ and $X_i^*$ is an odd cycle.
    \end{enumerate}
\end{lemma}

\begin{figure}[htb]
    \centering
    \includegraphics[height=9cm]{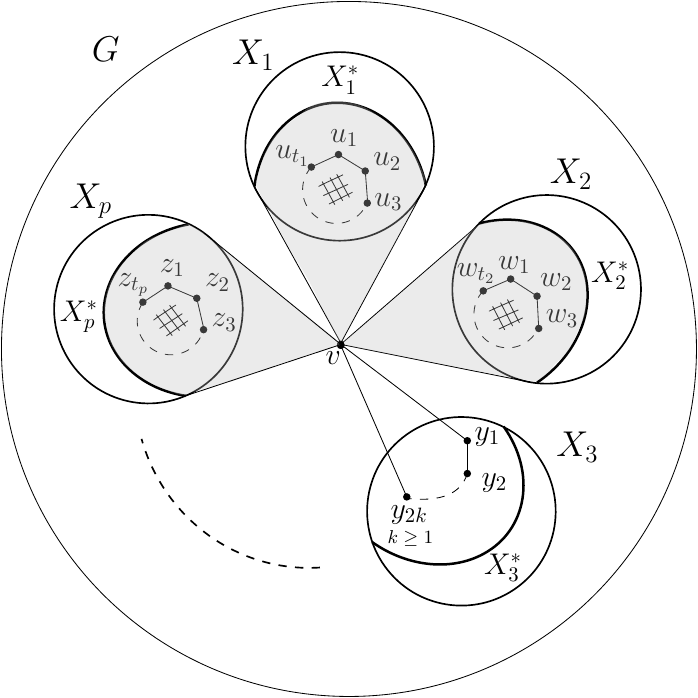}
    \caption{An example of a proper $(t_1, \dots, t_p)$-partition where all $t_i=\chi(X_i)\geq 3$ except for $t_3 = \chi(X_3) = 2$; note that each of the components $X_i^*$ are cliques, except for $X_3^*$, which is an odd cycle.}
    \label{general proper Mozhan}
  \end{figure}

In terms of how to picture proper $(t_1, \ldots, t_p)$-partitions, see Figure \ref{general proper Mozhan} for an example where each of the components $X_i^*$ are cliques, except for $X_3^*$, which is an odd cycle.

If a vertex-critical graph $G$ has $\chi\geq 10$, then we can always get a $(t_1, t_2, t_3)$-partition of $G$ where $t_i\geq 3$ for all $i$ (via our above discussion) and hence Lemma \ref{main mozhan lemma}(3) always gives a clique instead of odd cycle. This is why it is convenient to assume $\chi\geq 10$ when using the method of Mozhan partitions, but we will not have that luxury in this paper.

In addition to Lemma \ref{main mozhan lemma}, we will also use two other lemmas from \cite{DefsPaper} concerning these partitions.

\begin{lemma}
    [Kostochka, Rabern, and Stiebitz \cite{DefsPaper}]
    \label{shared neighbourhoods}
    Let $G$ be a vertex-critical graph and let $(v, X_1, \ldots, X_p)$ be a proper $(t_1, \dots, t_p)$-partition of $G$.
    If $y \in V(X_i^*)\setminus v$ is a low vertex of $G$ for some $i$, and $y$ has a neighbour in $X_j^*\setminus \{v\}$ for some $j\neq i$, then $N_{X_j}(v) = N_{X_j}(y)$.   
\end{lemma}

Let $P=(v, X_1, \dots , X_p$) be an optimal $(t_1, \dots , t_p)$-partition of $G$ and say vertex $y \in X_j$ for some $j$. Let $Y_i = X_i$ if $i \neq j$ and let $Y_j = (X_j \cup \{v\}) /\{y\}$; we say that $P' = (y, Y_1, \dots, Y_p)$ is obtained by \emph{swapping} $v$ with $y$ in $P$. We can see that $P'$ is a $(t_1, \dots, t_p)$-partition of $G$ if and only if $\chi(G[Y_j]) = t_j$. The following lemma (in fact a combination of two lemmas from \cite{DefsPaper}) provides conditions for when swapping maintains the property of being optimal, or even being proper.

\begin{lemma} [Kostochka, Rabern, and Stiebitz \cite{DefsPaper}]
\label{swapping is okay!}
    Let $G$ be a vertex-critical graph, let $P = (v, X_1, \dots , X_p)$ be a $(t_1, \dots, t_p)$-partition of $G$, and let $y\in V(X_j^*)\setminus v$. 

\begin{enumerate}
    \item If $P$ is optimal and $d_{X_j}(v) = t_j$, then swapping $v$ for $y$ in $P$ yields an optimal $(t_1, \dots, t_p)$-partition of $G$.
    \item If $P$ is proper and $y$ is low, then swapping $v$ for $y$ in $P$ yields a proper $(t_1, \dots, t_p)$-partition of $G$. 
\end{enumerate}

\end{lemma}

\section{Proofs of Theorems \ref{fullstrong} and \ref{minus4strong}}

We now present the proof of Theorem \ref{fullstrong}, which we restate here for convenience.

\setcounter{theorem}{2}
\begin{theorem} 
    Let $G$ be a vertex-critical graph with $\chi= \Delta = 9$ which does not contain any of $F, R_1, R_4, R_5, R_7, R_8, R_9$. Then $G$ contains $K_3 \lor E_6$ as a subgraph.
\end{theorem}

\begin{proof} Since $G$ is vertex-critical with $\chi=9$, we may choose an optimal $(2,3,3)$-partition of $G$, say $P = (v, X, Y, Z)$. 

\begin{claim}
    \label{special vertex is low}
   We may assume that the special vertex $v$ of $P$ is low, and hence that $P$ is proper. 
\end{claim}

\begin{proofc} By Lemma \ref{main mozhan lemma}(1), $v$ has at least two edges into $X$ and at least three edges into both $Y$ and $Z$.  Suppose, for a contradiction, that $v$ is high -- then it has degree 9, so it has one additional edge to some part $X, Y$, or $Z$. Since $x$ does not meet the hypothesis for Lemma \ref{main mozhan lemma}(3) for some part, the corresponding graph may not be a clique or an odd cycle. However, Lemma \ref{main mozhan lemma}(3) does apply for the other two parts and hence the corresponding graphs there are cliques or odd cycles. But then $G$ contains the basic structure of $F$ with $v_3=v$ a high vertex. Since $G$ does not contain $F$, at least one of $v_1, v_2, v_4, v_5, v_6$ must be a low vertex.  But then by Lemma~\ref{swapping is okay!}(1), we can swap $v$ with this low neighbour, obtaining an optimal (2,3,3)-partition where the special vertex is low.
\end{proofc}

Consider now Algorithm 1. This algorithm performs a sequence of vertex swaps to gradually change our proper (2,3,3)-partition $P=(v, X, Y, Z)$, while maintaining that we indeed have a proper (2,3,3)-partition. Once a pair of vertices are swapped in the algorithm, we say that both of them have \emph{moved}, a designation they will carry for the remainder of the algorithm. Each time a vertex moves, we create a new partition. We begin the algorithm with our original $P$ and denote the partition when a certain vertex $v_i$ is special by $P_i=(v_i, X_{v_i}, Y_{v_i}, Z_{v_i})$. Furthermore, we use $X_{v_i}^*, Y_{v_i}^*, Z_{v_i}^*$ to denote the components of $G[X_{v_i} \cup \{v_i\}], G[Y_{v_i} \cup \{v_i\}], G[Z_{v_i} \cup \{v_i\}]$ containing $v_i$, respectively. Note that in our original partition $P=(v, X, Y, Z)$ (which we label as $P=(v_1, X_{v_1}, Y_{v_1}, Z_{v_1})$), since $v=v_1$ is low, Lemma~\ref{main mozhan lemma}(2) and Lemma~\ref{main mozhan lemma}(3) implies that $X_{v_1}^*$ is an odd cycle and that $Y_{v_1}^*, Z_{v_1}^*$ are cliques (each of size four). The algorithm only ever swaps a pair of low vertices so by Lemma \ref{swapping is okay!}(2), $(v_i, X_{v_i}, Y_{v_i}, Z_{v_i})$ will be a proper (2,3,3)-partition throughout the algorithm, and hence by Lemma \ref{main mozhan lemma}(3) we know $X_{v_i}^*$ will be an odd cycle and $Y_{v_i}^*, Z_{v_i}^*$ will be cliques (each of size four) throughout the algorithm. Given any $P_i$, we may refer to $X_{v_i}$ as the \emph{$X$-set} or the \emph{cycle part}, while we may refer to $Y_{v_i}, Z_{v_i}$ as the the $Y$-set or $Z$-set, respectively, or as the \emph{clique parts}. 

Algorithm 1 begins by swapping the special vertex $v=v_1$ with any low neighbour of $v_1$. After that, the algorithm alternates between a swap involving the $Y$-set, and a swap involving either the $X$-set or the $Z$-set (with a preference for the $X$-set). Note that the condition for swapping with an $X$-set or $Z$-set is stricter than the condition for swapping with the $Y$-set. See Figure \ref{Algorithm 1} for an example of of a sequence of swaps in Algorithm 1.

\begin{tcolorbox}
 \textbf{Algorithm 1.}

 \vspace*{.1in}

 $\bullet$ \emph{Initialize:} $P_1 = P$, $v_1=v$, and $i=1$.

 \vspace*{.1in}

 $\bullet$ \emph{Iteration 1:} In $P_1$, swap $v_1$ with a low neighbour $w_1$ in one of $X_{v_1}, Y_{v_1}, Z_{v_1}$, setting $j_1 =1, 2,$ or $3$, respectively. Set $v_2= w_1$ and let $P_2$ be the resulting proper (2,3,3)-partition $(v_2, X_{v_2}, Y_{v_2}, Z_{v_2})$.

\vspace*{.1in}

$\bullet$ \emph{Iteration $i$, $i\geq 2$}: 
\begin{enumerate}
\item[(1)] \emph{If $j_{i-1}\in\{1, 3\}$:} \\
In $P_i$, if there exists a low vertex $w_i \in Y_{v_i}^*\setminus \{v_i\}$  that has never moved, then swap $v_i$ with $w_i$, set $j_i=2$, set $v_{i+1}=w_i$, let $P_{i+1}$ be the resulting proper (2,3,3)-partition, and iterate. Otherwise, terminate the algorithm.
\item[(2)] \emph{If $j_{i-1}=2$:}\\
In $P_i$, if every low neighbour of $v_i$ in $X_{v_i}^*$ has never moved, then swap $v_i$ with some such low neighbour $w_i$,  set $j_i=1$, set $v_{i+1}=w_i$, let $P_{i+1}$ be the resulting proper (2,3,3)-partition, and iterate. Otherwise, if every low vertex in $Z_{v_i}^*\setminus \{v_i\}$ has never moved, then swap $v_i$ with some such low vertex $w_i$, set $j_i=3$, set $v_{i+1}=w_i$, let $P_{i+1}$ be the resulting proper (2,3,3)-partition, and iterate. Otherwise, terminate the algorithm.
\end{enumerate}
\end{tcolorbox}

\begin{figure}[htb]
    \centering
   \includegraphics[height=4.5cm]{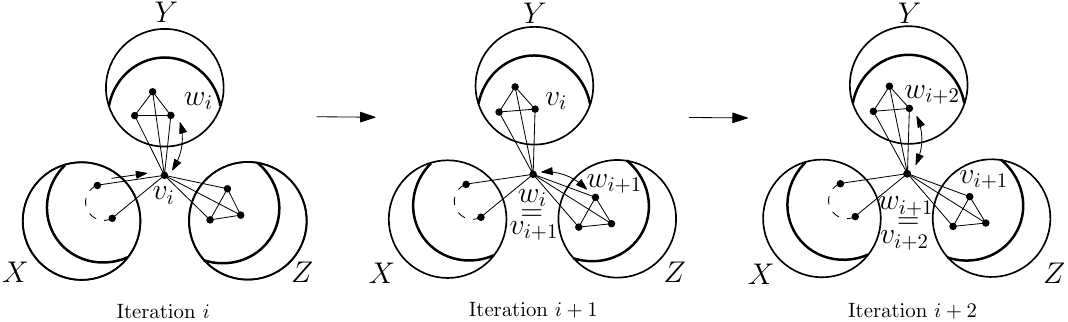}
    \caption{Example of three successive iterations of Algorithm 1: in iteration $i$ we imagine that $v_i$ has just been swapped out of the $X$-set in the previous step, so we swap $v_i$ with some $w_i\in Y$; in iteration $i+1$ we swap the new $v_{i+1}$ with $w_{i+1}\in Z$ (meaning $v_{i+1}$ must have a low neighbour in the $X$-set that has already moved), and; in iteration $i+2$ we swap the new $v_{i+2}$ with $w_{i+2}\in Y$.}
    \label{Algorithm 1}
  \end{figure}

Algorithm 1 makes sense as written because we have forbidden $R_1, R_5$ and $R_7$. First, we need to ensure that in iteration 1, the vertex $v_1$ \emph{has} a low neighbour, and forbidding $R_1$ accomplishes this. Secondly, in iteration $i(2)$, we need to ensure that the condition on the low vertices in the $X$-set or $Z$-set is not met vacuously, that is, we need to ensure that if one of these conditions is met then there is at least one low neighbour in the $X$-set or $Z$-set, respectively, to swap with. So if the special vertex $v_i$ has just swapped out of the $Y$-set, we want it to have at least one low neighbour in both the $X$-set and the $Z$-set. In this scenario, $v_i$ has a low neighbour in the $Y$-set (the vertex it just swapped with), and forbidding $R_5$ and $R_7$ forces $v_i$ ($v_6$) to have a low neighbour in both the $X$-set and the $Z$-set.

In Algorithm 1, each time we swap a pair of vertices, we are swapping the special vertex with a vertex that has never moved before. Since the graph is finite, this means that the algorithm must terminate eventually.  We divide the remainder of our proof into two cases, depending on whether our algorithm terminated due to the stopping condition in (1) or in (2) above. In both cases we will show the existence of the desired subgraph $K_3 \lor E_6$. \\

\noindent\textbf{Case 1:} \textit{At termination, $j_{i-1}=2$.}

In this case, at termination, we know that $v_i$ was swapped out of the $Y$-set in the previous step (since $j_{i-1}=2$). Moreover, we know that the special vertex $v_i$ has a low neighbour $a\in X_{v_i}$ and a low neighbour $b \in Z_{v_i}$, both of which have already moved. Suppose that $a$ moved before $b$. (The argument of this case is almost exactly the same if $b$ moved before $a$, so we omit those details to avoid repetition.) Since $a$ moved before $b$, we know that $b$ cannot be the initial special vertex $v_1$. So given the algorithm's alternating pattern, and the fact that $b \in Z_{v_i}$, we know that $b$ must have swapped out of the $Y$-set in a previous step.  Let $b'$ be the vertex that was special when this swap occurred; note that after this swap $b'$  remains in the $Y$-set until the algorithm terminates.

Consider the point in the algorithm when $a$ was special. Since $a$ moves before both $v_i$ and $b$, it must be that $v_i, b \in Y_a$. Furthermore, since $a \sim v_i \sim b$, we have $v_i, b \in Y_a^*$. In addition, $a$ must have one more neighbour in the clique $Y_a^*$, call it $m$. So $V(Y_a^*)=\{a, v_i, b, m \}$.

\begin{claim} \label{m does not move after a}
    The vertex $m$ does not move after $a$ moves.
\end{claim}

\begin{proofc}
    Assume for contradiction that $m$ did move at some point in the algorithm after $a$ moves. 
    First, we assume that $m$ moves after $a$ but before $b$. Then consider the point in the algorithm when $m$ is special (so immediately after it swapped out of the $Y$-set). Since $m \sim a$ and $a$ moves into the $X$-set, $a \in X_m^*$. Since $a$ is a neighbour of $m$ that has already moved, $m$ must move into the $Z$-set. But then when $b$ becomes the special vertex later in the algorithm, we would have $a \in X_b^*$ and $m \in Z_b^*$ (since $b \sim a, m$), and both $a,m$ have moved. So the algorithm should have terminated when $b$ was the special vertex, contradiction. We may now assume that $m$ moves after both $a$ and $b$. But then when $m$ is special, we would already have $a \in X_m^*$ and $b \in Z_m^*$, and the algorithm should have terminated at this point, contradiction.
\end{proofc}

\begin{claim} \label{v_i and b}
    $v_i \sim b'$.
\end{claim}

\begin{proofc}
    Consider the point in the algorithm when $b$ was special. We know that $v_i$ has not moved out of the $Y$-set yet (since it is special at termination). So since $b \sim v_i$, we get $v_i \in Y_b^*$. In addition, $b'$ has just swapped into the $Y$-set (because $b'$ is the vertex that swaps with $b$ to make $b$ special). So $b' \in Y_b^*$ as well, meaning that $v_i \sim b'$.
\end{proofc}

\begin{claim} \label{mb'}
    $m, b' \in Y_{v_i}^*$.
\end{claim}
\begin{proofc}
Consider the final partition $P_i$ at termination. Since $m \sim v_i$ and $m$ does not move after $a$ (by Claim \ref{m does not move after a}), we get that $m \in Y_{v_i}^*$. On the other hand, since $v_i \sim b'$ (by Claim \ref{v_i and b}) and since $b'$ remains in the $Y$-set until termination, we also get that $b' \in Y_{v_i}^*$. 
\end{proofc}

Let us now assemble everything we know about the partition $P_i$ with special vertex $v_i$. We know by assumption that $a\in X_{v_i}^*$, $b \in Z_{v_i}^*$, and Claim \ref{mb'} tells us that $m, b' \in Y_{v_i}^*$. Let us label now the other vertices adjacent to $v_i$ in $X_{v_i}^*, Y_{v_i}^*$ and $Z_{v_i}^*$: let $u_1$ be this other vertex in the $X$-set, let $u_2$ be this other vertex in the $Y$-set, and let $u_3, u_4$ be these other two vertices in the $Z$-set.  See the leftmost picture in Figure \ref{getting final adjacencies (1)}, where edges we currently know exist are indicated with thin edges, and edges we are about to find (in Claim \ref{claim for conclusion adjacencies}) are indicated with thick edges. The thin edges in the picture come from $X_{v_i}^*, Y_{v_i}^*, Z_{v_i}^*$ as just described, plus the triangle $(a, b, m)$ which is a part of the clique $Y_a^*$, and the edge $b\sim b'$. Our final goal is the $K_3\vee E_6$ pictured on the right of Figure  \ref{getting final adjacencies (1)}. In particular, we have already established that $\{v_i, m, b\}$ induces a $K_3$, and we have already established that $v_i$ is adjacent to all of the vertices in our proposed $E_6$ (namely $a, b', u_1, u_2, u_3, u_4$). 
We already know that $m\sim a, b', u_2$ and that $b\sim a, u_3, u_4$. So it remains only for us to show that $m \sim u_1, u_3, u_4$ and that $b\sim u_1, u_2$. That is, Claim \ref{claim for conclusion adjacencies} will complete this case (of $j_{i-1}=2$).

\begin{figure}[htb]
    \centering
   \includegraphics[height=5.2cm]{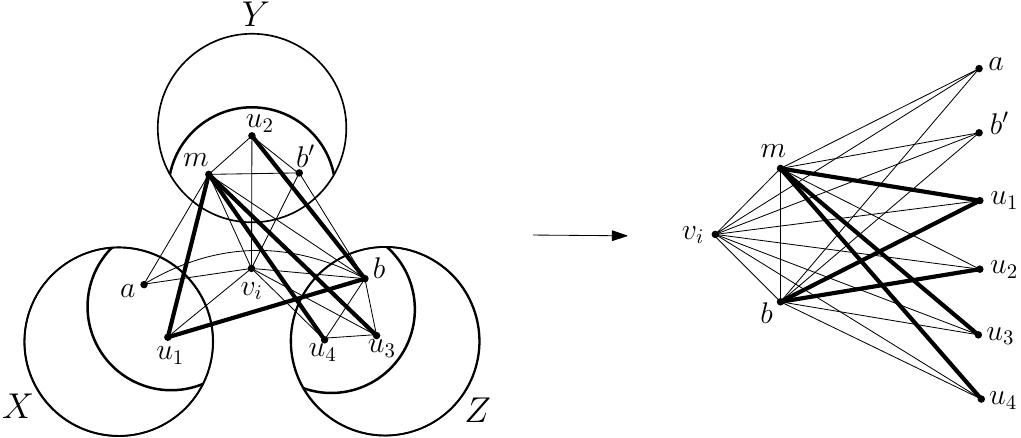}
    \caption{On the left we see the partition $P_i$ at termination; the thick edges are those found in Claim \ref{claim for conclusion adjacencies}, while the thin edges are those found previously. The graph to the right illustrates the final result.}
    \label{getting final adjacencies (1)}
  \end{figure}

\begin{claim}
    \label{claim for conclusion adjacencies}
$m\sim u_1, u_3, u_4$ and $b\sim u_1, u_2$. 
\end{claim}

\begin{proofc} In order to show that $b \sim u_2$, we apply Lemma~\ref{shared neighbourhoods}: since $b$ is low and $b\in Z_{v_i}^*\setminus \{v_i$\}, and since $b\sim m$  with $m \in Y_{v_i}^*\setminus \{v_i\}$, we have that $N_{Y_{v_i}}(v_i) = N_{Y_{v_i}}(b)$, and thus $b\sim u_2$ since $v_i\sim u_2$ (see Figure \ref{getting final adjacencies (1)}). 

In order to show $m \sim u_3,u_4$,  we again apply Lemma~\ref{shared neighbourhoods}: $m$ is low and $m\in Y_{v_i}^*\setminus \{v_i\}$, and since $m\sim b$ (by the last paragraph) and $b \in Z_{v_i}^* \setminus \{v_i\}$, we have that $N_{Z_{v_i}}(v_i) = N_{Z_{v_i}}(m)$, and thus $m\sim u_3,u_4$ since $v_i\sim u_3, u_4$ (see Figure \ref{getting final adjacencies (1)}). 

It remains now to show that $u_1 \sim m, b$, for which we again apply Lemma~\ref{shared neighbourhoods}: $m,b$ are both low with $m\in Y_{v_i}^*\setminus \{v_i\}$ and $b\in Z_{v_i}^*\setminus \{v_i\}$, and since $m,b\sim a$  with $a\in X_{v_i}^* \setminus \{v_i\}$, we have that $N_{X_{v_i}}(v_i) = N_{X_{v_i}}(m)=N_{X_{v_i}}(b)$, and thus $m,b\sim u_1$ since $v_i\sim u_1$ (see Figure \ref{getting final adjacencies (1)}). 
\end{proofc}

\vspace*{.1in}

\noindent\textbf{Case 2:} \textit{At termination, $j_{i-1}\in\{1, 3\}$.}

In this case, at termination, we know that $v_i$ was either swapped out of the $X$-set or the $Z$-set in the previous step (since $j_{i-1} \in \{1,3\}$). Moreover, every low neighbour of $v_i$ in $Y_{v_i}$ has already moved. Since $v_i$ either swapped out of the $X$-set or the $Z$-set, it certainly has a low neighbour in whichever set it just swapped out of (the vertex that it swapped with to become special). Since $G$ does not contain $R_4, R_5$, we cannot have that all the neighbours of $v_i$ in the $Y$-set are high.
Since $G$ does not contain $R_8, R_9$, we cannot have that exactly two of the neighbours of $v_i$ in the $Y$-set are high. Thus, $v_i$ has two low neighbours in $Y_{v_i}$, call them  $a', b'$. By the assumption on this case, $a',b'$ have already moved; let $a,b$ be the vertices that $a',b'$ swapped with, respectively, when they swapped into the $Y$-set. Note that $a, b\neq v_i$, since by the assumption of this case, $v_i$ just swapped out of the $X$-set or the $Z$-set, not the $Y$-set. Without loss, say that $a'$ swapped into the $Y$-set before $b'$. In fact, assume that $a'$ was the first of any of the low vertices in $Y^*_{v_i}$ to swap into the $Y$-set, followed by $b'$ second. In particular, this means that $a'$ was already in the $Y$-set when $b'$ moved in.

\begin{claim}
    \label{thm 1 case 2 adjacencies 1}
    The vertices $\{a, a', b, b'\}$ induce a clique with the possible exception of the edge $ab'$. (In fact we will find later that $ab'$ exists.)
\end{claim}

\begin{proofc}
Edges $aa'$ and $bb'$ exist by definition of $a, b$. As $a', b' \in Y_{v_i}^*$, which is a clique, they are adjacent. Consider the point at which $b'$ was special: at this point we must have had $a', b$ in the $Y$-set, and since both are adjacent to $b'$, we get that $a', b\in Y^*_{b'}$ which is a clique, so $a'\sim b$.  Similarly, we consider the point where $a'$ was special: at this point we must have had that $a, b$ in the $Y$-set, and since both are adjacent to $a'$, we get that $a, b\in Y^*_{a'}$ which is a clique, so $a\sim b$. 
\end{proofc}

We know that $a,b$ were swapped out of the $Y$-set by $a',b'$, respectively. The following claim gives more information about what happened immediately after each one of those swaps. 

\begin{claim}
\label{a into X claim}
    After $a$ swaps out of the $Y$-set, it must swap into the $X$-set. After $b$ swaps out of the $Y$-set, it must swap into the $Z$-set. 
\end{claim}

\begin{proofc} When $a$ is special, it has a low neighbour $a'$ in the $Y$-set. Since $G$ does not contain $R_7$, it must be that $a$ has some low neighbour in the $X$-set at this point. Since $a$ just swapped out of the $Y$-set (it swapped with $a'$), the algorithm will look to move $a$ into the $X$-set if possible. If it is not possible, it is because some low neighbour of $a$ in the $X$-set has already moved; assume for a contradiction that this is the case. Since we guaranteed that this condition cannot be met vacuously, we may assume that when $a$ was special, it had some low neighbour $c$ in $X^*_a$ which had already moved, and so $a$ swapped into the $Z$-set.

Consider the point in the algorithm when $c$ was special. This must be before $a, b$ swap out of the $Y$-set. At this point, $a \in Y_{c}^*$, and as $a$ and $b$ are adjacent, $b \in Y_{c}^*$, so we have $b \sim c$. Consider now the (later) point in the algorithm when $b$ was special. This was after $a, c$ had both moved. At this point, $c \in X_b^*$ and $a \in Z_b^*$. Then since $a,c$ are low neighbours of $b$ that have already moved, the algorithm should have terminated here by the stopping condition in (2). This is a contradiction since we assumed the algorithm terminated by the stopping condition in (1). So the first sentence of our claim statement is true.

Consider now the point in the algorithm when $b$ was special. We know that $b$ just swapped out of the $Y$-set (swapped with $b'$), so the algorithm will first look to swap $b$ into the $X$-set if possible. But this is not possible here, because we know that $a$ is in the $X$-set at this point (by the first sentence of our claim), and $a\sim b$ by the previous claim, and $a$ has already moved. So, the algorithm will next look to swap $b$ into the $Z$-set if possible. But we know it must be possible, since the algorithm doesn't terminate at this point (it terminates when $v_i\neq b$ is the special vertex). So we also get that the second sentence of our claim statement is true. 
\end{proofc}

Consider the point in the algorithm when $b'$ is special, just before it swaps $b$ out of the $Y$-set. We know that $a'$ must have swapped into the $Y$-set at some prior step (swapping $a$ out). 
Since $b'\sim a', b$ by Claim \ref{thm 1 case 2 adjacencies 1}, we get that $a', b \in Y_{b'}^*$. Then $Y_{b'}^*$ is a $K_4$ made of $b'$, $a', b$, and some other vertex $m$. We now analyze this vertex $m$.

\begin{claim}
    \label{m never moves}
       The vertex $m$ never moves, and $m\sim a$.
\end{claim}

\begin{proofc}
Assume to the contrary that $m$ did move at some step before termination. First we assume $m$ has already moved when $b'$ is special. Then since $m \in Y_{b'}$, it must remain in the $Y$-set until the algorithm terminates. At termination, $a', b' \in Y_{v_i}^*$ and since $a',b' \sim m,$ we have $m \in Y_{v_i}^*$. But since $m$ swapped into the $Y$-set before $b'$, we contradict our assumption that $b'$ was the second low neighbour of $v_i$ in $Y_{v_i}$ to swap into the $Y$-set, after $a'$. 

We next assume for contradiction that $m$ moved after $b'$ was special. We consider the partition when $a'$ was special; we know this is before $b'$ was special because we have assumed that $a$ leaves the $Y$-set before $b$. At this step, none of the vertices $a, b, m$ have moved yet, so they are all in the $Y$-set. Then $a, b, m \in Y_{a'}^*$, so in particular $m\sim a,b$. Since we are assuming that $m$ moves after $b'$, it must move after $a$ and $b$, so we consider the later point in the algorithm when $m$ is special. By Claim \ref{a into X claim}, $a$ would be in the $X$-set and $b$ in the $Z$-set by this point. Since $m$ just moved out of the $Y$-set, the algorithm will look to swap $m$ into either the $X$-set or $Z$-set. But since $m$ has a low neighbour in both of these sets that has already moved ($a,b$), the algorithm would terminate by the stopping condition in (2), which is a contradiction since we assumed the algorithm terminated by the stopping condition in (1).

Since we have established that $m$ never moves, consider the point at which $a'$ is special. Then $a\in Y^*_{a'}$, and since $m\sim a'$ and $m\in Y_{a'}$, we get that $m\in Y^*_{a'}$, and hence $m\sim a$.
\end{proofc}

We now return our focus to the final partition at termination. Here $v_i$ is special with $a', b' \in Y_{v_i}^*$ and $a, b,$ in the $X$-set and $Z$-set, respectively, by Claim \ref{a into X claim}. Since $m$ never moves, it remains in the $Y$-set, and since $m \sim a', b'$, we have $m \in Y_{v_i}^*$. In order to achieve our final structure, we propose making one additional swap. We swap $v_i$ with $a'$, making $a'$ the special vertex. As $a'$ is low, this new partition is still proper by Lemma \ref{swapping is okay!}(2). Call this new partition $(a', X_{a'}, Y_{a'}, Z_{a'})$.

Let us now assemble everything we know about this new partition with special vertex $a'$. We know that $v_i, m, b' \in Y_{a'}^*$, and we know that $a \in X_{a'}^*$ since $a \sim a'$ by definition and $b \in Z_{a'}^*$ since $a' \sim b$ by Claim \ref{thm 1 case 2 adjacencies 1}. Let us label now the other vertices adjacent to $a'$ in $X_{v_i}^* $ and $ Z_{v_i}^*$: let $u_1$ be this other vertex in the $X$-set and let $u_2, u_3$ be these other two vertices in the $Z$-set.  See the leftmost picture in Figure \ref{getting final adjacencies (2)}, where edges we currently know exist are indicated with thin edges, and edges we are about to find (in Claim \ref{thm 1 case 2 final adjacencies}) are indicated with thick edges. The thin edges in the picture come from $X_{v_i}^*, Y_{v_i}^*, Z_{v_i}^*$ as just described, plus $m \sim a$ (by Claim \ref{m never moves}),  $a \sim b$ (by Claim \ref{thm 1 case 2 adjacencies 1}), and $b \sim b'$. Our final goal is the $K_3 \vee E_6$ pictured on the right of Figure \ref{getting final adjacencies (2)}. In particular, we have already established that $\{a', b', b\}$ induces a $K_3$ by Claim \ref{thm 1 case 2 adjacencies 1}, and we have already established that $a'$ is adjacent to all of the vertices in our proposed $E_6$ (namely $v_i, m, a, u_1, u_2, u_3$). 
We already know that $b' \sim v_i, m, b$ and that $b \sim a, u_2, u_3, b'$. So it remains only for us to show that $b' \sim a, u_1, u_2, u_3$ and that $b\sim v_i, m, u_1$. That is, Claim \ref{thm 1 case 2 final adjacencies} will complete this case (of $j_{i-1} \in \{1,3\}$).

\begin{figure}[htb]
    \centering
   \includegraphics[height=5.2cm]{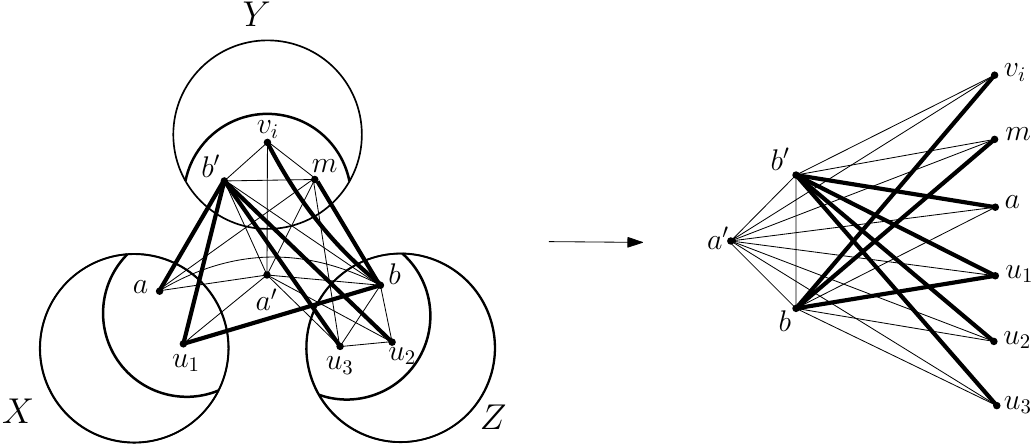}
    \caption{On the left we see the partition $(a', X_{a'}, Y_{a'}, Z_{a'})$; the thick edges are those found in Claim \ref{thm 1 case 2 final adjacencies}, while the thin edges are those found previously. The graph to the right illustrates the final result. }
    \label{getting final adjacencies (2)}
  \end{figure}

\begin{claim} \label{thm 1 case 2 final adjacencies}
    $b' \sim a, u_1, u_2, u_3$ and $b \sim v_i, m, u_1$.
\end{claim}

\begin{proofc}
     In order to show that $b \sim v_i, m$, we apply Lemma \ref{shared neighbourhoods}: since $b$ is low and $b \sim b'$ with $b' \in Y_{a'}^*\setminus \{a'\}$, we have that $N_{Y_{a'}}(a') = N_{Y_{a'}}(b)$, and thus $b \sim v_i, m$ because $a' \sim v_i, m$. (See Figure \ref{getting final adjacencies (2)}).

In order to show that $b' \sim u_2, u_3$, we apply Lemma \ref{shared neighbourhoods} again: since $b'$ is low and $b'\sim b$ with $b \in Z_{a'}^*\setminus \{a'\}$, we have that $N_{Z_{a'}}(a') = N_{Z_{a'}}(b')$, and thus $b' \sim u_2, u_3$ because $a' \sim u_2, u_3$. (See Figure \ref{getting final adjacencies (2)}).

To show that $b' \sim a$, we again apply Lemma \ref{shared neighbourhoods}: since $a$ is low and $a\sim m$ with $m \in Y_{a'}^*\setminus \{a'\}$, we have that $N_{Y_{a'}}(a') = N_{Y_{a'}}(a)$, and thus $a\sim b'$ since $a' \sim b'$. (See Figure \ref{getting final adjacencies (2)}).

It remains to show that $b, b' \sim u_1$, for which we again apply Lemma \ref{shared neighbourhoods}: since $b, b'$ are both low and $b, b'\sim a$ (by our last line) with $a \in X_{a'}^*\setminus \{a'\}$, we have that $N_{X_{a'}}(a') = N_{X_{a'}}(b') = N_{X_{a'}}(b)$, and thus $b, b'\sim u_1$ since $a' \sim u_1$. (See Figure \ref{getting final adjacencies (2)}).
\end{proofc}
\end{proof}

We now present the proof of Theorem \ref{minus4strong}, which is restated below. 

\setcounter{theorem}{4}
\begin{theorem}
    Let $G$ be a vertex-critical graph with $\chi= \Delta = 9$ which does not contain any of $F, R_2, R_5$. Then $G$ contains $(K_3 \lor E_6)^{-4}$ as a subgraph.
\end{theorem}

\begin{proof}
Since $G$ is vertex-critical with $\chi = 9$, we may choose an optimal (2,3,3)-partition of $G$, say $P= (v, X, Y, Z)$. We use the same argument as in Claim~\ref{special vertex is low} of Theorem~\ref{fullstrong} to assume that the special vertex $v$ is low and hence that $P$ is a proper $(2,3,3)$-partition of $G$. 

Consider now Algorithm 2, which we will discuss using the same notation and terminology as Algorithm 1. In fact, the two algorithms are very similar, with two key differences. The first difference is that switches into the $X$-set (i.e. the cycle part) are not allowed in Algorithm 2. In particular, this affects iteration 1, where we look to swap into either the $Y$-set or the $Z$-set, and also iteration i(2), where we look only to swap into the $Z$-set. The second key difference between Algorithm 1 and Algorithm 2 is that in iteration i(1), we raise the bar for the swapping condition into the $Y$-set -- we now need \emph{every} low vertex in $Y_{v_i^*}\setminus v_i$ to have never moved, rather than just having one that has never moved. Note that this means we have symmetric swapping conditions in $i(1)$ and $i(2)$, which is different from Algorithm 1.

Algorithm 2 makes sense as written because we have forbidden $R_2$ and $R_5$. First, we need to ensure that in iteration 1, the vertex $v_1$ \emph{has} a low neighbour in one of the clique parts, and forbidding $R_2$ accomplishes this. Secondly, in iteration $i$, we need to ensure that the condition on the low vertices in the $Y$-set or the $Z$-set is not met vacuously, that is, we need to ensure that if the condition in (1) or (2) is met, then there is at least one low vertex in the $Y$-set or the $Z$-set, respectively, to swap with. So if the special vertex $v_i$ has just swapped out of the $Y$-set, we want it to have at least one low neighbour in the $Z$-set. In this scenario, $v_i$ has a low neighbour in the $Y$-set (the vertex it just swapped with), and forbidding $R_5$ forces $v_i$ $(v_6)$ to have a low neighbour in the $Z$-set. Similarly, a special vertex $v_i$ which has just swapped out of the $Z$-set will be forced to have a low neighbour in the $Y$-set.

In Algorithm 2, just as in Algorithm 1, we get termination because swaps must always involve a vertex that has never moved before. In the proof of Theorem \ref{fullstrong}, we needed two different cases to find our target structure, depending on whether Algorithm 1 stopped due to condition (1) or condition (2). However now in Algorithm 2, these two conditions are symmetric, and indeed the sets $Y, Z$ can be relabelled without loss of generality if we wish. Hence, we may assume without loss that Algorithm 2 terminates by the stopping condition in (1), that is, when $j_{i-1} = 2$.

\begin{tcolorbox}
 \textbf{Algorithm 2.}

 \vspace*{.1in}

 $\bullet$ \emph{Initialize:} $P_1 = P$, $v_1=v$, and $i=1$.

 \vspace*{.1in}

$\bullet$ \emph{Iteration 1:} In $P_1$, swap $v_1$ with a low neighbour $w_1$ in $Y_{v_1}$ or $Z_{v_1}$, setting $j_1 =1$ or $2$, respectively. Set $v_2 = w_1$ and let $P_2$ be the resulting proper $(2,3,3)$-partition $(v_2, X_{v_2}, Y_{v_2}, Z_{v_2})$.

 \vspace*{.1in}

$\bullet$ \emph{Iteration $i$, $i \geq 2$}:

\begin{enumerate}
\item[(1)] \emph{If $j_{i-1} = 2$:} 
\\
In $P_i$, if every low vertex in $Y_{v_i}^* \setminus v_i$  has never moved, then swap $v_i$ with some such low vertex $w_i$, set $j_i=1$, set $v_{i+1}=w_i$, let $P_{i+1}$ be the resulting proper partition $(2,3,3)$-partition, and iterate. Otherwise, terminate the algorithm.

\item[(2)] \emph{If $j_{i-1}=1$:}
\\
In $P_i$, if every low vertex in $Z_{v_i}^* \setminus v_i$ has never moved, then swap $v_i$ with some such low vertex $w_i$, set $j_i=2$, set $v_{i+1}=w_i$, let $P_{i+1}$ be the resulting proper partition $(2,3,3)$-partition, and iterate. Otherwise, terminate the algorithm.

\end{enumerate}
\end{tcolorbox}

Consider the special vertex $v_i$ at termination. Since $j_{i-1}=2$, $v_i$ has just swapped out of the $Z$-set and has a low neighbour $a' \in Y_{v_i}^*$ which has already moved.  Let one of the other vertices in $Y_{v_i}^* \setminus \{a'\}$ be called $m$. Then, in particular, $m \sim a'$. When $a'$ moved into the $Y$-set in some prior step, it must have swapped a vertex out, call it $a$. Then because of our algorithm's alternating pattern, we know $a$ must have swapped into the $Z$-set, where it remains at termination.

\begin{claim}
    \label{thm 2 w doesn't move, and edge am}
    The vertex $m$ never moves, and $a \sim m$.
\end{claim}

\begin{proofc}

Assume to the contrary that $m$ moves at some step in the algorithm. We know that $m$ is in the $Y$-set at termination, so it must swap into the $Y$-set either before or after $a'$ does. 

First assume for a contradiction that $m$ swapped into the $Y$-set before $a'$ did. We consider the step in the algorithm when $a'$ was special, just before it swapped into the $Y$-set. Since $m \sim a'$, we have that $m \in Y_{a'}^*$, but then the algorithm should have terminated at this step by the stopping condition in (1) since $m$ has already moved. Since $a'\neq v_i$, this is a contradiction. 

Now assume for a contradiction that $m$ swapped into the $Y$-set after $a'$ did. We consider the step in the algorithm when $m$ was special, just before it swapped into the $Y$-set. We find the same contradiction as in the first case: $m \sim a'$ implies that $a' \in Y_m^*$. Then since $a'$ has already moved, the algorithm should have terminated at this point, but $m \neq v_i$, contradiction. Thus $m$ never moves.

Since $m$ never moves, it resides in the $Y$-set for every step of the algorithm. Thus, it is in the $Y$-set with $a$ before $a$ moves. In particular, since $a\sim a'\sim m$, we get $a, m \in Y_{a'}^*$, so $a \sim m$. 
\end{proofc}

We consider again the final partition at termination. Here $v_i$ is special, $m, a' \in Y_{v_i}^*$, and $a$ is in the $Z$-set. In order to achieve our final structure, we perform one additional swap. We swap $a'$ with $v_i$ so that $a'$ is special and $v_i$ is in the $Y$-set, and $a,m$ remain where they are. Since $a'$ is low, this new partition is still proper by Lemma \ref{swapping is okay!}(2). Call this new partition $(a', X_{a'}, Y_{a'}, Z_{a'})$.

We now assemble everything we know about this new partition with special vertex $a'$. We know that  $v_i, m \in Y_{a'}^*$ and $a \in Z_{a'}^*$. Let us label now the other vertices adjacent to $a'$ in $X_{a'}^*, Y_{a'}^*$ and $Z_{a'}^*$: let $u_1, u_2$ be these other vertices in the $X$-set, let $u_3$ be this other vertex in the $Y$-set, and let $u_4, u_5$ be these other two vertices in the $Z$-set. As $Y_{a'}^*$ and $Z_{a'}^*$ are cliques, we have that $\{m, v_i, u_3\}$ and $\{a, u_4, u_5\}$ induce cliques, and we also know that $a \sim m$ by Claim \ref{thm 2 w doesn't move, and edge am}. Our final goal is a $(K_3 \vee E_6)^{-4}$. In particular, we will show that $\{a', a, v_i\}$ can play the role of $K_3$ for this structure, and that $\{m, u_1, \dots, u_5\}$ can play the role of $E_6$, with the four potentially missing edges being those between $\{u_1, u_2\}$ and $\{v_i, a\}$. We have already established that $a'$ is adjacent to all of the vertices in our proposed $E_6$ (namely $m, u_1, \dots, u_5)$. We also already know that $v_i \sim m, u_3$ and that $a \sim m, u_4, u_5$. So it remains only for us to show that $a \sim v_i, u_3$ and $v_i \sim u_4, u_5$. That is, Claim \ref{minus4 final adjacencies} will complete the proof of Theorem \ref{minus4strong}.

\begin{claim} \label{minus4 final adjacencies}
    $a \sim v_i, u_3$ and $v_i \sim u_4, u_5$.
\end{claim}

\begin{proofc}
 In order to show that $a \sim v_i, u_3$, we apply Lemma \ref{shared neighbourhoods}: since $a $ is low and  $a \sim m$ with $m \in Y_{a'}^*\setminus \{a'\}$, we have that $N_{Y_{a'}}(a') = N_{Y_{a'}}(a)$, and thus $a \sim v_i, u_3$ since $a' \sim v_i, u_3$. To show that $v_i \sim u_4, u_5$, we apply Lemma \ref{shared neighbourhoods} a second time: since $v_i$ is low and  $v_i \sim a$ (by our previous line) with $a \in Z_{a'}^* \setminus \{a'\}$, we have that $N_{Z_{a'}}(a') = N_{Z_{a'}}(v_i)$, and so $v_i \sim u_4, u_5$ because $a' \sim u_4, u_5$.
\end{proofc}
\end{proof}

\section{Proof of Theorem \ref{minus6strong}}

We now present the proof of Theorem \ref{minus6strong}, which we restate here for convenience.

\setcounter{theorem}{5} 
\begin{theorem}
    Let $G$ be a vertex-critical graph with $\chi= \Delta = 9$ which does not contain any of $F, R_1, R_3, R_{10}$. Then $G$ contains $(K_3 \lor E_6)^{-6}$ as a subgraph.
\end{theorem}

\begin{proof} 
Since $G$ is vertex-critical with $\chi = 9$, we may choose an optimal (2,3,3)-partition of $G$, say $P= (v, X, Y, Z)$. We use the same argument as in Claim~\ref{special vertex is low} of Theorem~\ref{fullstrong} to assume that the special vertex $v$ is low and hence that $P$ is a proper $(2,3,3)$-partition of $G$. 

Consider now Algorithm 3, which we will discuss using the same notation and terminology as Algorithms 1 and 2. Algorithm 1 allows the special vertex to swap with any of the three parts ($X$-set, $Y$-set, $Z$-set) in the first iteration, and then alternates between swaps with the $Y$-set and swaps with the $X$-set or $Z$-set. Algorithm 2 maintains the same pattern, except it allows no swaps at all with the $X$-set. Algorithm 3 again allows swaps with any of the three parts, but it strictly prioritizes swaps with the clique sets ($Y$-set and $Z$-set), preferring to alternate swaps between these two sets only, and if forced to swap with the $X$-set, looking to terminate the algorithm immediately unless an ``ideal next swap'' presents itself. Let us now say more about how this all works.

\begin{tcolorbox}
 \textbf{Algorithm 3.}

 \vspace*{.1in}

 $\bullet$ \emph{Initialize:} $P_1 = P$, $v_1=v$, and $i=1$.

 \vspace*{.1in}

$\bullet$ \emph{Iteration 1:}
 In $P_1$, swap $v_1$ with a low neighbour $w_1$ in one of $X_{v_1}, Y_{v_1}, Z_{v_1}$, prioritizing a $w_1$ in a clique part. Moreover, if $w_1$ is from a clique part, prioritizing $w_1$ which has a low neighbour in the \emph{other} clique part. Set $j_1 = 1, 2,$ or $3$, respectively, set $v_2 = w_1$ and let $P_2$ be the resulting proper (2,3,3)-partition $(v_2, X_{v_2}, Y_{v_2}, Z_{v_2})$. 

 \vspace*{.1in}

$\bullet$ \emph{Iteration $i$, $i\geq 2$}:
\begin{enumerate}
\item[(1)] \emph{If $j_{i-1} = 2$ $[j_{i-1} = 3]$}:
\
\begin{enumerate}
    \item  In $P_i$, if there is a low vertex in $Z_{v_i}^* \setminus v_i $ $[Y_{v_i}^* \setminus v_i]$: 
    
    If every low vertex in $Z_{v_i}^* \setminus v_i $ $[Y_{v_i}^* \setminus v_i]$ has never moved, then swap $v_i$ with such a low $w_i$, prioritizing a $w_i$ which has a low neighbour in $Y_{v_i}[Z_{v_i}]$. Set $j_i = 3(2),$ set $v_{i+1} = w_i$, let $P_{i+1}$ be the resulting proper (2,3,3)-partition, and iterate. Otherwise, terminate the algorithm.

    \item In $P_i$, if there is no low vertex in $Z_{v_i}^* \setminus v_i $ $[Y_{v_i}^* \setminus v_i]$:
    
    If every low neighbour of $v_i$ in $X_{v_i}^*$ has never moved, then swap $v_i$ with such a low $w_i$, set $j_i =1$, $v_{i+1} = w_i$, let $P_{i+1}$ be the resulting proper (2,3,3)-partition, and iterate. Otherwise, terminate the algorithm.
\end{enumerate}

\item[(2)] \emph{If $j_{i-1} = 1$:} 
\\
    In $P_i$, if every low vertex in $(Y_{v_i}^* \cup Z_{v_i}^*) \setminus v_i$ has never moved, then swap $v_i$ with such a low $w_i$, prioritizing a $w_i$ which has a low neighbour in the opposite clique part. Set $j_i = 2$ or $3$ if $w_i$ is in $Y_{v_i}$ or $Z_{v_i}$, respectively. Set $v_{i+1} = w_i$, let $P_{i+1}$ be the resulting proper (2,3,3)-partition, and iterate. Otherwise, terminate the algorithm.

\end{enumerate}
\end{tcolorbox}

In iteration 1 of Algorithm 3, $v_1$ swaps into the $X$-set only if it has no low neighbours in either the $Y$-set or $Z$-set (where Algorithm 1 expressed no preference between the three sets). Moreover, even if Algorithm 3 is able to swap into the Y-set or Z-set, it prioritizes specially picking a vertex to swap with that will allow the next iteration to go smoothly as well: namely, if possible it chooses a vertex in the $Y$-set or $Z$-set that has a low neighbour in the \emph{other} clique part. Note that this means in iteration 2, the special vertex will be able to swap with this neighbour in the other clique part, beginning this desired alternating pattern of swaps with the $Y$-set and $Z$-set.

Suppose we are at a point in Algorithm 3 where we have just been forced to swap into the $X$-set (i.e. $j_{i-1}=1$), then we set the bar very high for continuing the algorithm: we only swap (as opposed to terminate) if every low vertex in \emph{all} of $(Y^*_{v_i}\cup Z^*_{v_i})\setminus v_i$ has never moved. If this bar is met, then just as in iteration 1, we choose (if possible) a vertex in the $Y$-set or $Z$-set that has a low neighbour in the \emph{other} clique part.

In any iteration $i\geq 2$ where $j_{i-1}\neq 1$ (i.e. we have not just swapped with the $X$-set), we have the same condition for swapping into a set. In particular, whether a vertex is going to swap into the $X$-set, $Y$-set, or $Z$-set, it must satisfy the following condition: no low neighbour of $v_i$ in $X^*_{v_i}$, $Y^*_{v_i}$, or $Z^*_{v_i}$, respectively can have already moved. This is actually the condition for swaps that Algorithm 2 has; in Algorithm 1 there was a lower bar to meet for swapping into the $Y$-set. Also, similarly to iteration 1, if $v_i$ is going to swap into the $Y$-set or $Z$-set, then Algorithm 3 prioritizes specially picking a vertex to swap with that has a low neighbour in the \emph{other} clique part.

Algorithm 3 makes sense as written because we have forbidden $R_1, R_3, R_{10}$. First, we need to ensure that in iteration 1, the vertex $v_1$ \emph{has} a low neighbour, and forbidding $R_1$ accomplishes this. 
Secondly, in iteration $i(1)(b)$, we need to ensure that the condition on the low vertices in the $X$-set is not met vacuously, i.e. that there is some low vertex in the $X$-set to swap with. The situation here is that $v_i$ has just swapped out of either the $Y$-set or the $Z$-set (so it has one low neighbour there) but has no low neighbours in the other clique part. So forbidding $R_{10}$ forces $v_i$ ($v_6$) to have a low neighbour in the $X$-set. Finally, for iteration $i(2)$, we'll ensure that a special vertex $v_i$ which has just swapped out of the $X$-set has at least one low neighbour in the $Y$-set or the $Z$-set. In this scenario, $v_i$ has a low neighbour in the $X$-set (the vertex it just swapped with), and forbidding $R_{3}$ forces $v_i$ ($v_6$) to have a low neighbour in either the $Y$-set or the $Z$-set.

In Algorithm 3, as in Algorithms 1 and 2, we get termination because swaps must always involve a vertex that has never moved before. We divide the remainder of our proof into two cases. If the algorithm terminates due to the stopping condition in $(1b)$, we will find a contradiction. If it terminates due to the stopping condition in either $(1a)$ or $(2)$, we show that $G$ contains $(K_3 \lor E_6)^{-6}$ as desired. \\

\noindent\textbf{Case 1}: \textit{Algorithm 3 terminates due to the stopping condition in (1b).}

In this case, at termination, we know that $v_i$ was either just swapped out of the $Y$-set or the $Z$-set (since $j_{i-1} \in \{2,3\}$), and that $v_i$ has a low neighbour $a \in X_{v_i}$ that has already moved. Given the pattern of Algorithm 3, we know that $a$ either swapped out of one of the clique parts in a previous step, or that $a$ was the initial special vertex $v_1$. We will consider both possibilities, and in fact show that neither is possible.

We first assume, for a contradiction, that $a$ was in one of the clique parts before swapping into the $X$-set.  We may assume, without loss of generality, that $a$ was in the $Y$-set previously. We claim that this means $v_i$ was not in the $Y$-set previously. If it was, then $a,v_i$ were both in the $Y$-set before $a$ moved. Let $a'$ be the vertex that swapped into the $Y$-set and swapped $a$ out. Consider the point in the algorithm when $a'$ was special. At this point, $a \in Y_{a'}^*$ and $a \sim v_i$, so $v_i \in Y_{a'}^*$, and thus $v_i \sim a'$. Then after $a'$ swapped into the $Y$-set (swapped with $a$), we know that $a'$ remained in the $Y$-set until termination. Now we look to a later iteration. We consider the vertex that swapped into the $Y$-set and swapped $v_i$ out, call this vertex $v_i'$. At the step when $v_i'$ was special, $v_i \in Y_{v_i'}^*$, and since $v_i \sim a'$, we also have $a' \in Y_{v_i'}^*$, so $v_i' \sim a'$. But since $a'$ is a low vertex that has already moved, the algorithm should have terminated at this point by the stopping condition in either $(1a)$ or (2), contradicting our assumption that it stopped by the condition in $(1b)$. So we establish our claim, and hence we know that $v_i$ was in the $Z$-set prior to becoming the special vertex (not the $Y$-set). But then consider the step in the algorithm when $a$ was special, just after $a$ swapped out of the $Y$-set. At this point, $v_i \in Z_a^*$ since $v_i \sim a$ and $a$ moves before $v_i$. Since $a$ has a low neighbour ($v_i$) in the $Z$-set, the algorithm should have either swapped $a$ with a low neighbour in the $Z$-set or terminated, as described in $(1a)$. Since we know that $a$ actually swapped into the $X$-set, we have a contradiction.

We now assume, for a contradiction, that $a$ was the initial special vertex $v_1$. Since $j_{i-1} \in \{2, 3\}$, we know that $v_i$ was in one of the clique parts before becoming special; in particular, $v_i$ was in one of the clique parts in the initial partition $P_1 = (a=v_1, X_{v_1}, Y_{v_1}, Z_{v_1})$. So then when $a$ was special, it had a low neighbour $v_i$ in one of the clique parts. Thus, iteration 1 dictates that $a$ should have swapped into either the $Y$-set or the $Z$-set. But we know that $a \in X_{v_i}^*$ at termination, so $a$ swapped into the $X$-set, contradiction.\\

\noindent\textbf{Case 2}: \textit{Algorithm 3} terminates due to the stopping condition in $(1a)$ or $(2).$ 

In this case, at termination, $v_i$ has just swapped out of either the $X$-set, $Y$-set, or the $Z$-set (since $j_{i-1} \in \{1, 2, 3\}$). Regardless, by either of these stopping conditions, $v_i$ has a low neighbour $a'$ in either the $Y$-set or the $Z$-set that has already moved. Without loss of generality suppose that $a' \in Y_{v_i}^*$.  Then since the algorithm stops with $v_i$ looking into the $Y$-set at $a'$, we know that $v_i$ cannot have swapped out of the $Y$-set, so $j_{i-1} \in \{1, 3\}$ and $v_i$ either swapped out of the $X$-set or the $Z$-set. Let one of the other neighbours of $v_i$ in $Y_{v_i}^* \setminus \{a'\}$ be called $m$. Let the vertex that $a'$ swapped with when it moved into the $Y$-set be called $a.$ Note that $a$ moves before $v_i$ in the course of the algorithm.
 
 \begin{claim}
    \label{thm 3 case 2 adjacencies 1}
     The vertices $\{m, a', v_i\}$ induce a clique. Also, $a \sim m, $ $a \sim a'$, and $m$ never moves.
 \end{claim}
 
 \begin{proofc}
     
We know that these three vertices induce a clique because $m, a', v_i \in Y_{v_i}^*$. Also $a' \sim a$ by definition. Applying the same argument as in Claim~\ref{thm 2 w doesn't move, and edge am} of the proof of Theorem~\ref{minus4strong} shows that $a\sim m$ and that $m$ never moves. 
 \end{proofc}

Since $a'$ swapped into the $Y$-set and swapped $a$ out, we know that $a$ either swapped into the $X$-set or the $Z$-set and remained there until termination. The argument is easier in the latter case, where we in fact find a slightly stronger subgraph.

\begin{claim}If $a$ is in the $Z$-set at termination, then $G$ contains $(K_3 \lor E_6)^{-4}$. \end{claim}

\begin{proofc}
We consider the partition $(v_i, X_{v_i}, Y_{v_i}, Z_{v_i})$ at termination. Here we know $a', m \in Y_{v_i}^*$ (by definition) and $a \in Z_{v_i}$ by assumption. To achieve our final structure, we propose making one additional swap: we swap $v_i$ with $a'$, making $a'$ special and swapping $v_i$ into the $Y$-set. Note that after this swap, $m$ and $a$ are still in the $Y$-set and $Z$-set respectively. Since $a'$ is low, this new partition is proper by Lemma \ref{swapping is okay!} (2). Call the new partition $(a', X_{a'}, Y_{a'}, Z_{a'})$.  

We now assemble everything that we know about the partition $(a', X_{a'}, Y_{a'}, Z_{a'})$. By Claim \ref{thm 3 case 2 adjacencies 1}, $m, v_i \in Y_{a'}^*$ and $a \in Z_{a'}^*$. We now label the other vertices adjacent to $a'$ in $X_{a'}^*, Y_{a'}^*,$ and $Z_{a'}^*$: let $u_1, u_2$ be these other vertices in the $X$-set, let $u_3$ be this other vertex in the $Y$-set, and let $u_4, u_5$ be these other vertices in the $Z$-set. As $Y_{a'}^*$ and $Z_{a'}^*$ are cliques, we have that $\{v_i, m, u_3\}$ and $\{a, u_4, u_5\}$ induce cliques. We also have $a \sim m$ by Claim \ref{thm 3 case 2 adjacencies 1}. 

We will show that $\{a, a', v_i\}$ induce a $K_3$ and that $\{m, u_1, \dots u_5\}$ can play the role of $E_6$, with the four potentially missing edges being those between $\{v_i, a\}$ and $\{u_1, u_2\}$. We have already established that $a'$ is adjacent to all of the vertices in our proposed $E_6$ (namely $m, u_1, \dots u_5$). We already know that $v_i \sim m, u_3$ and $a \sim u_4, u_5, m$. It remains for us to show that $a \sim v_i, u_3$ and $v_i \sim u_4, u_5$. 

In order to show that $a \sim v_i, u_3$, we apply Lemma \ref{shared neighbourhoods}: since $a$ is low and $a \sim m$ with $m \in Y_{a'}^*\setminus \{a'\}$, we have that $N_{Y_{a'}}(a') = N_{Y_{a'}}(a)$, and thus $a \sim v_i, u_3$ since $a' \sim v_i, u_3$. To show that $v_i \sim u_4, u_5$, we apply Lemma \ref{shared neighbourhoods} again: since $v_i$ is low and $v_i \sim a$ (by our previous line) with $a \in Z_{a'}^* \setminus \{a'\}$, we have that $N_{Z_{a'}}(a') = N_{Z_{a'}}(v_i)$, and so $v_i \sim u_4, u_5$ because $a' \sim u_4, u_5$.
\end{proofc}

We may now assume that $a$ is in the $X$-set at termination. We will show that $G$ contains  $(K_3 \lor E_6)^{-6}$ via a series of claims.

\begin{claim}
    \label{thm 3 subcase 1 adjacency 1}
    $a \sim v_i$.
\end{claim}

\begin{proofc}
    In $(v_i, X_{v_i}, Y_{v_i}, Z_{v_i})$, we propose swapping $v_i$ with $a'$ so that $a'$ is special and $v_i$ is in the $Y$-set. This swap creates a new proper (2,3,3)-partition because $a'$ is low (Lemma \ref{swapping is okay!}(2)). Note that after this swap, we still have $a$ in the $X$-set and $m$ in the $Y$-set. Call this new partition $(a', X_{a'}, Y_{a'}, Z_{a'})$. Since $a' \sim a, m, v_i$ by Claim \ref{thm 3 case 2 adjacencies 1}, we have $a \in X_{a'}^*$ and $m, v_i \in Y_{a'}^*$. Applying Lemma \ref{shared neighbourhoods}, since $a$ is low and $a \sim m$ (by Claim \ref{thm 3 case 2 adjacencies 1}) with $m \in Y_{a'}^*$, we have that $N_{Y_{a'}}(a') = N_{Y_{a'}}(a)$, and so $a \sim v_i$ because $a' \sim v_i$.
\end{proofc}

\begin{claim}
    Before $v_i$ moved, it was in the $X$-set.
    \label{x is in X}
\end{claim}

\begin{proofc} We already know that $j_{i-1} \in \{1, 3\}$; assume for contradiction that $v_i$ was in the $Z$-set before moving. We consider the point in the algorithm when $a$ was special, just after it swapped out of the $Y$-set. At this point, $v_i \in Z_a^*$ (since $a \sim v_i$ by our prior claim and $a$ moves before $v_i$). Thus, the partition $(a, X_a, Y_a, Z_a)$ meets the condition in 1(a): $j_{i-1} = 2$ and there is a low vertex $v_i \in Z_{a}^* \setminus \{a\}$. So this iteration of the algorithm should have either ended with termination or with $a$ swapping into the $Z$-set, contradiction.
\end{proofc}

Now consider again the step where $a$ was the special vertex. We know that $v_i$ is in the $X$-set before it moves (Claim \ref{x is in X}) and that $a$ moves before $v_i$, so $v_i$ is in $X_a$, and since $a \sim v_i$ by Claim \ref{thm 3 subcase 1 adjacency 1}, we get that $v_i$ is a neighbour of $a$ in $X_a^*$. Let $u_1$ be the other neighbour of $a$ in the $X$-set.

\begin{claim} \label{a swaps with vi}
    When $a$ swapped into the $X$-set, it swapped with $v_i$.
\end{claim}

\begin{proofc}
    Assume, for a contradiction, that when $a$ moved into the $X$-set it swapped with $u_1$. See Figure \ref{path in X argument} for a depiction of both before and after this swap. Consider the $X$-set when $u_1$ is special (the right picture in Figure \ref{path in X argument}); recall that $X_{u_1}^*$ is an induced odd cycle. Let  $p$ be the other neighbour of $u_1$ in $X$ (besides $a$), and let $P$ be the induced path between $a$ and $p$ in the $X$-set. Note that $P$ contains $v_i$.
    
    As the algorithm progresses, $P$ will remain unchanged in the $X$-set until some vertex $s$ becomes special and wants to switch with one of the vertices in $P$. But since $X_s^*$ must be an induced cycle, this means that the two neighbours of $s$ in the $X$-set must be $a, p$. At this point however, the algorithm would have stopped according to 1(b), since $a$ has already moved. Hence $P$ remains unchanged in the $X$-set at termination. This is a contradiction, since $v_i$ is a member of $P$, and we know that $v_i$ is the special vertex at termination.
\end{proofc}

\begin{figure}[htb]
    \centering
   \includegraphics[height=2.5cm]{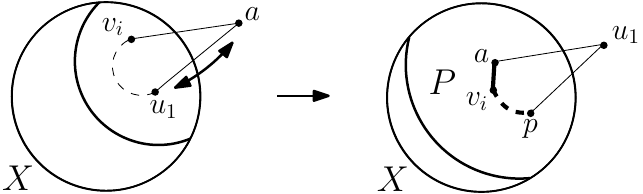}
    \caption{Two points in the algorithm as described in Claim \ref{a swaps with vi}.}
    \label{path in X argument}
  \end{figure}

\begin{claim}
    \label{thm 3 subcase 1 adjacencies 2}
    $v_i \sim u_1$.
\end{claim}

\begin{proofc}
After $a$ swaps with $v_i$ (Claim \ref{a swaps with vi}), $v_i$ becomes special, and the resulting partition is actually the partition at termination: $P_i = (v_i, X_{v_i}, Y_{v_i}, Z_{v_i})$. Here we know that $a, u_1$ are in the $X$-set  and $a' \in Y_{v_i}^*$ (by definition of $u_1, a'$). We know that $a\sim u_1$; to show that $v_i$ is adjacent to $u_1$, we suggest a sequence of swaps between low vertices, each of which creates a new proper (2,3,3)-partition by Lemma \ref{swapping is okay!} (2). 

First we swap $v_i$ into the $Y$-set and swap $a'$ out. Then since $a' \sim a$ by definition, we can swap $a'$ into the $X$-set and swap $a$ out. 
Call this new partition $(a, \tilde{X}_{a}, \tilde{Y}_a, \tilde{Z}_a)$; 
note we are using the tilde notation here so as not to confuse this new partition (which is created after termination) 
with the partition $(a, X_{a}, Y_a, Z_a)$ (the partition when $a$ was special during Algorithm 3). 

Then $ a', u_1 \in \tilde{X}_a^*$ (since $a \sim a', u_1$ by definition) and $v_i \in Y_a^*$ by Claim \ref{thm 3 subcase 1 adjacency 1}. Since $v_i$ is low and $v_i \sim a' \in X_a^* \setminus \{a\}$ by definition, we apply Lemma~\ref{shared neighbourhoods} to find that $N_{X_a}(a) = N_{X_a}(v_i)$, and thus $v_i \sim u_1$ since $a \sim u_1$. 
\end{proofc}

We now assemble everything we know about the partition $(a, X_a, Y_a, Z_a)$ when $a$ was special during Algorithm 3 (just after $a'$ swapped $a$ out of the $Y$-set). 
At this point, we know that two neighbours of $a$ in $X_a^*$ are $v_i, u_1$ (see discussion above Claim 16). 
We also know that $a' \in Y_a^*$ since $a'$ has just swapped into the $Y$-set and swapped $a$ out. Moreover, we know that $m \in Y_a^*$ since $a \sim m$ and $m$ never moves (Claim \ref{thm 3 case 2 adjacencies 1}).
We now label the other vertices adjacent to $a$ in $Y_a^*$ and $Z_a^*$: let $u_2$ be this other vertex in the $Y$-set and let $u_3, u_4, u_5$ be these other vertices in the $Z$-set. 

As $Y_a^*$ is a clique, we have that $\{a', m, u_2\}$ induce a clique. We also have that $a' \sim v_i$ by definition and $v_i \sim m, u_1$ by Claims \ref{thm 3 case 2 adjacencies 1} and \ref{thm 3 subcase 1 adjacencies 2}. Our final goal is a $(K_3 \vee E_6)^{-6}$. In particular, we have already established that $\{a, a', v_i\}$ induce a $K_3$, and we will show that $\{m, u_1, \dots u_5\}$ can play the role of $E_6$, with the six potentially missing edges being those between $\{v_i, a'\}$ and $\{u_3, u_4, u_5\}$. 
We have already established that $a$ is adjacent to all of the vertices in our proposed $E_6$ (namely $m, u_1, \dots u_5$). We already know that $v_i \sim m, u_1$ and $a' \sim m, u_2$. So it remains only for us to show that $v_i \sim u_2$ and $a' \sim u_1$. 

 In order to show that $v_i \sim u_2$, we apply Lemma \ref{shared neighbourhoods}: since $v_i$ is low and $v_i \sim a'$ with $a' \in Y_{a}^*\setminus \{a\}$, we have that $N_{Y_{a}}(a) = N_{Y_{a}}(v_i)$, and thus $v_i \sim u_2$ since $a \sim u_2$. To show that $a' \sim u_1$, we apply Lemma \ref{shared neighbourhoods} once more: since $a'$ is low and $a' \sim v_i$ with $v_i \in X_{a}^* \setminus \{a\}$, we have that $N_{X_{a}}(a) = N_{X_{a}}(a')$, and so $a' \sim u_1$ because $a \sim u_1$.
\end{proof}

\section{An infinite family}

In this final section we present an infinite family of graphs that is as we promised in the introduction: each member of the family will have $\chi=\Delta=9$ with average degree above $8.75$, but will contain none of the substructures listed in Table 1.

\begin{figure}[htb]
    \centering
   \includegraphics[height=10cm]{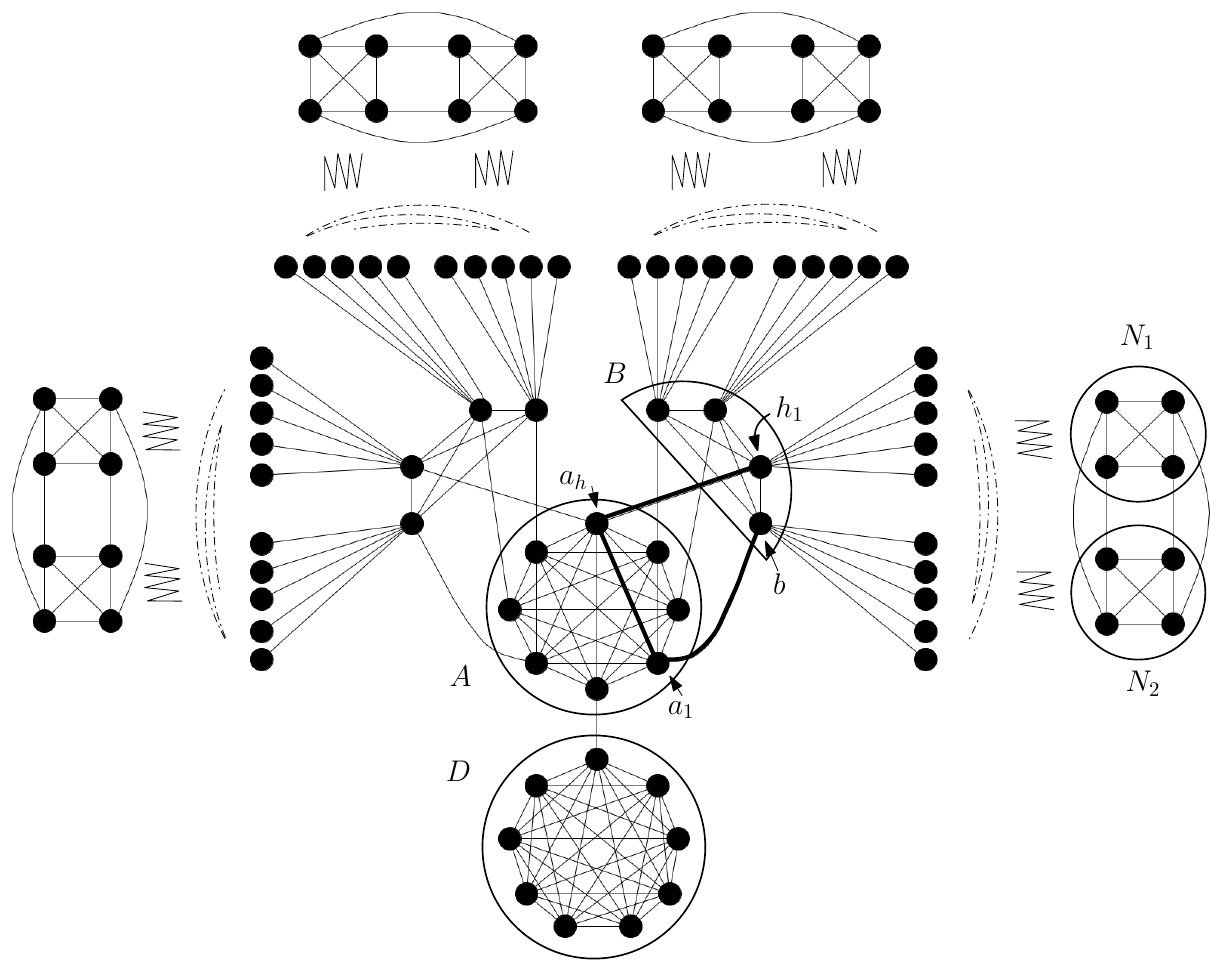}
    \caption{An $H^*$ building block. Note that a solid jagged line between two sets of vertices indicates a complete bipartite graph between those sets; a dashed jagged line indicates that the two sets are joined by a complete bipartite graph less one perfect matching. Note that the labelled vertices $b, h_1, a_1, a_h$ and the bolded segment they lie on will be used much later.}
    \label{example graph no connection}
  \end{figure}

Our infinite family of graphs will be called $\mathcal{H}$. Each member of $\mathcal{H}$ will be made using a basic building block $H^*$, which is the graph depicted in Figure \ref{example graph no connection}. In these figures, a solid jagged line between two sets of vertices indicates a complete bipartite graph between those sets; a dashed jagged line indicates that the two sets are joined by a complete bipartite graph less one perfect matching.

\begin{figure}[htb]
    \centering
   \includegraphics[height=4.2cm]{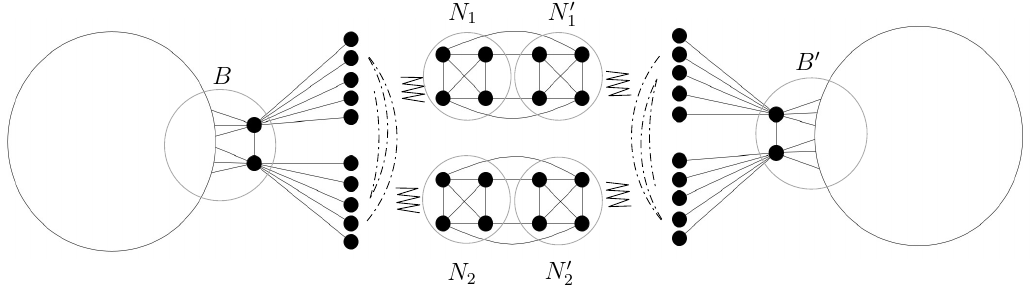}
    \caption{A connection point between two $H^*$ building blocks. The two blocks are connected by edges between $N_1$ and $N_1'$ and between $N_2$ and $N_2'$. The open circles on either side represent the rest of the $H^*$ blocks which are not pictured in detail.}
    \label{example graph with connection}
  \end{figure}

To see how we can use copies of $H^*$ to make members of our family $\mathcal{H}$, consider the four pairs of $K_4$'s occurring along the outside boundary of $H^*$ in Figure \ref{example graph no connection}; we will consider each such pair of $K_4$'s as ``free for connection''. 
Take any such free pair of $K_4$'s (say with vertex sets named $N_1, N_2$) and delete the 4 edges joining them. Then take another copy of $H^*$ and similarly delete the edges from one of its free $K_4$ pairs (say with vertex sets named $N_1', N_2'$). Then we connect these two $H^*$ building blocks by joining $N_1, N_1'$ as a pair (with four matching edges, as used to connect $N_1$ to $N_2$), and joining $N_2, N_2'$ as a pair (with four matching edges, as used to connect $N_1'$ to $N_2'$). See Figure \ref{example graph with connection}. Note that after this procedure, we no longer consider the vertex sets $N_1, N_1', N_2, N_2'$ as corresponding free $K_4$ pairs, but the resulting graph has exactly 6 pairs of $K_4$'s on the boundary of its structure, and we consider those pairs free for connections. Our family $\mathcal{H}$ consists of $H^*$ plus  any number of copies of $H^*$ that have been iteratively joined via the above procedure.

Consider any $H\in\mathcal{H}$. Note that $\chi(H)=9$ since it contains a $K_9$ (on the vertex set $D$ of $H^*$, indicated in Figure  \ref{example graph no connection}). All the vertices in $H^*$ have either degree 9 or 8, with the only low vertices being: the vertices of $A$ (with the exception of its topmost vertex in Figure \ref{example graph no connection}), and; the vertices of $D$ (with the exception of its topmost vertex in Figure \ref{example graph no connection}). Hence the average degree of $H^*$ is just above 8.845 ($H^*$ contains 15 low vertices and 82 high). The degrees of all vertices remain the same after our combination procedure above, so this average degree is the same for $H$. (If one desired an average degree that was even higher, it is possible to use a second type of building block in order to make this happen, but we have not included those details here for the sake of brevity.)

Note that $H^*$ is certainly \emph{not} vertex-critical, since one of the vertices in $D\subseteq H^*$ is a cut-vertex, and deleting this vertex from $H^*$ creates two components, both of which have chromatic number $8$ (note both contain $K_8$). On the other hand, if we perform this same deletion in any $H\in \mathcal{H}\setminus H^*$, then at least one of the two components created does contain a $K_9$ and hence the resulting graph still has chromatic number 9. 

It remains to show that if $H\in \mathcal{H}$, then $H$ contains none of our forbidden substructures. It is worth noting how close this is to failing: if we instead had two edges coming out of the $K_9$ in $H^*$, then the existence of the two high vertices in $D$ would mean we could find any of $R_7, R_8, R_9$ or $S$ within $D$.

\begin{lemma} \label{H doesnt have Rs}
   If $H\in \mathcal{H}$, then $H$ does not contain any of $R_1, \dots, R_{10}, Q,$ or $S$.
\end{lemma}

\begin{proof}
The structures $R_1, \dots, R_{10}, Q, S$ all contain an induced odd cycle where one of the vertices in the cycle is low, and its two neighbours on the cycle are high, so it suffices to show that $H$ contains no such structure. 

The only low vertices in $H$ come from the cliques $A$ and $D$ (in some copy of $H^*$). In $D$, the low vertices each have exactly one high neighbour (the top vertex of $D$ in Figure \ref{example graph no connection}), so it is not possible for any of these to be part of an induced odd cycle where their \emph{two} neighbours on the cycle are high. 

Let $a_1$ be a low vertex from $A$. Note that $a_1$ has exactly two high neighbours: the high vertex in $A$, which we'll call $a_h$, and; a neighbour outside of $A$, which we'll call $b$. Note that we may assume that $b\in B$ (see Figure \ref{example graph no connection}), since if $b\in D$, then $a_1, b$ are joined along a cut edge, which is certainly not part of any cycle in $H$.

Assume, for a contradiction, that there is an induced odd cycle $C$ in $H$ containing $a_h, a_1, b$ consecutively in $C$.  Note that since $A$ is a cut-set in the graph $H$, the other vertex of $C$ adjacent to $a_h$ (i.e besides $a_1$) must also be in $B$, say $h_1$. If not, then $C$ would have to include another vertex from $A$ (a clique) at some point, and this would contradict $C$ being an induced cycle. This means, without loss of generality, that our four consecutive vertices on $C$ ($h_1, a_h, a_1, b$) must be as labelled in Figure \ref{example graph no connection}. But then since $h_1, b$ are both members of $B$, they are adjacent, and hence we get a copy of $C_4$ induced by these four vertices. But since $C$ is supposed to be an induced odd cycle, this is a contradiction.
\end{proof}

It remains now to show that no $H\in \mathcal{H}$ contains the structure $F$ (a $K_4$ and an induced odd cycle meeting at a vertex, with all vertices high). There are a number of $K_4's$ in $H^*$ (and in any given $H\in \mathcal{H}$) but the all-high requirement rules out $K_4$'s from either $A$ or $D$ (see Figure \ref{example graph no connection}), and the attached induced odd cycle rules out a $K_4$ from $B$ (see Figure \ref{example graph no connection}). Therefore, if $H\in\mathcal{H}$ contains $F$, then the $K_4$-part of this $F$ comes from one of the $K_4$ pairs that was originally ``free for connection'' in some copy of $H^*$. If the $K_4$ in question was never \emph{used} for a connection, then it does not have an attached induced odd cycle (see Figure \ref{example graph no connection}). But this is also true if the $K_4$ \emph{was} used for a connection: see Figure \ref{example graph no connection}, where we can see that $N_1$ has no induced odd cycle attached to it.\\

\bibliographystyle{amsplain}
\bibliography{main}

\newpage

\end{document}